\documentclass[11pt]{article}
\usepackage{amssymb,amsfonts,amsmath,latexsym,epsf,tikz,url}

\newtheorem{theorem}{Theorem}[section]
\newtheorem{proposition}[theorem]{Proposition}

\newtheorem{corollary}[theorem]{Corollary}
\newtheorem{lemma}[theorem]{Lemma}
\newtheorem{remark}[theorem]{Remark}

\newcommand{\proof}{\noindent{\bf Proof.\ }}
\newcommand{\qed}{\hfill $\square$\medskip}

\begin{document}

\title{Distinguishing number and distinguishing index of  Kronecker product of two graphs}

\author{
Saeid Alikhani  $^{}$\footnote{Corresponding author}
\and
Samaneh Soltani
}

\date{\today}

\maketitle

\begin{center}
Department of Mathematics, Yazd University, 89195-741, Yazd, Iran\\
{\tt alikhani@yazd.ac.ir, s.soltani1979@gmail.com}
\end{center}


\begin{abstract}
The distinguishing number (index) $D(G)$ ($D'(G)$)  of a graph $G$ is the least integer $d$
such that $G$ has an vertex labeling (edge labeling)  with $d$ labels  that is preserved only by a trivial automorphism. The Kronecker product $G\times H$ of two graphs $G$ and $H$ is the graph with vertex set $V (G)\times
V (H)$ and edge set $\{\{(u, x), (v, y)\} | \{u, v\} \in E(G) ~and ~\{x, y\} \in  E(H)\}$.
In this paper we study the distinguishing number and the distinguishing index of Kronecker product of two graphs. 

\end{abstract}


\medskip
\noindent{\bf AMS Subj.\ Class.:} 05C15, 05C60. 

\section{Introduction and definitions}
Let $G=(V,E)$ be a simple graph of order $n\geqslant 2$. We use the the following notations: The set of vertices adjacent in $G$ to a vertex of a vertex subset $W\subseteq  V$ is the open neighborhood $N_G(W )$ of $W$. 
${\rm Aut}(G)$ denotes the automorphism group of $G$.  
A labeling of $G$, $\phi : V \rightarrow \{1, 2, \ldots , r\}$, is said to be $r$-distinguishing, 
if no non-trivial  automorphism of $G$ preserves all of the vertex labels.
The point of the labels on the vertices is to destroy the symmetries of the
graph, that is, to make the automorphism group of the labeled graph trivial.
Formally, $\phi$ is $r$-distinguishing if for every non-trivial $\sigma \in {\rm Aut}(G)$, there
exists $x$ in $V$ such that $\phi(x) \neq \phi(x\sigma)$. The distinguishing number of a graph $G$ is defined  by
\begin{equation*}
D(G) = {\rm min}\{r \vert ~ G ~\textsl{\rm{has a labeling that is $r$-distinguishing}}\}.
\end{equation*} 

This number has defined by Albertson and Collins \cite{Albert}. Similar to this definition, Kalinowski and Pil\'sniak \cite{Kali1} have defined the distinguishing index $D'(G)$ of $G$ which is  the least integer $d$
such that $G$ has an edge colouring   with $d$ colours that is preserved only by a trivial
automorphism. If a graph has no nontrivial automorphisms, its distinguishing number is  $1$. In other words, $D(G) = 1$ for the asymmetric graphs.
The other extreme, $D(G) = \vert V(G) \vert$, occurs if and only if $G = K_n$. The distinguishing index of some examples of graphs was exhibited in \cite{Kali1}. For 
instance, $D(P_n) = D'(P_n)=2$ for every $n\geqslant 3$, and 
$D(C_n) = D'(C_n)=3$ for $n =3,4,5$,  $D(C_n) = D'(C_n)=2$ for $n \geqslant 6$. It is easy to see that the value $|D(G)-D'(G)|$ can be large. For example $D'(K_{p,p})=2$ and $D(K_{p,p})=p+1$, for $p\geq 4$.  
A graph and its complement, always have the same automorphism group while their graph structure usually differs, hence $D(G) = D(\overline{G})$ for every simple graph $G$.
The distinguishing number and the distinguishing index of some graph products has been studied in literature (see \cite{Alikhani,Distlexico,Distneigh,Imrich,Imrich and  Klavzar}). 
The Cartesian product of graphs $G$ and $H$ is a graph, denoted $G\Box H$, whose vertex
set is $V (G) \times V (H)$. Two vertices $(g, h)$ and $(g', h')$ are adjacent if either $g = g'$ and
$hh' \in E(H)$, or $gg' \in E(G)$ and $h = h'$. Denote $G\Box G$ by $G^2$, and recursively define the
$k$-th Cartesian power of $G$ as $G^k = G\Box G^{k-1}$.
A non-trivial graph $G$ is called prime if $G = G_1\Box G_2$ implies that either $G_1$ or $G_2$ is
$K_1$. Two graphs
$G$ and $H$ are called relatively prime if $K_1$ is the only common factor of $G$ and $H$. We need knowledge of the structure of the automorphism group of the Cartesian product, which was determined by Imrich \cite{Imrich1969}, and independently by Miller \cite{Miller}.

\begin{theorem}{\rm \cite{Imrich1969,Miller}}\label{autoCartesian} Suppose $\psi$  is an automorphism of a connected graph $G$ with prime
	factor decomposition $G = G_1\Box G_2\Box \ldots \Box G_r$. Then there is a permutation $\pi$ of the set
	$\{1, 2, \ldots ,  r\}$ and there are isomorphisms  $\psi_i : G_{\pi(i)} \rightarrow G_i$, $i = 1,\ldots , r$, such that
	$$\psi(x_1, x_2 , \ldots , x_r) = ( \psi_1(x_{\pi(1)}),  \psi_2(x_{\pi(2)}), \ldots ,  \psi_r(x_{\pi(r)})).$$ 
\end{theorem}

The Kronecker product is one of the (four) most important graph products and seems to have been first introduced by K. \v Culik, who called it the cardinal product \cite{culik}. Weichsel \cite{Weichsel} proved that the Kronecker product of two nontrivial graphs is connected if and only if both factors are connected and at least one of them possesses an odd cycle. If both factors are connected and bipartite, then their Kronecker product consists of two connected components. 
The Kronecker product $G\times H$ of graphs $G$ and $H$ is the graph with vertex set $V (G)\times
V (H)$ and edge set $\{\{(u, x), (v, y)\} | \{u, v\} \in E(G) ~and ~\{x, y\} \in  E(H)\}$. The terminology
is justified by the fact that the adjacency matrix of a Kronecker graph product is given
by the Kronecker matrix product of the adjacency matrices of the factor graphs; see \cite{Weichsel} for details. However, this product is also known under several different names including categorical product, tensor product, direct product, weak direct product, cardinal product
and graph conjunction. The Kronecker product is commutative and associative in an obvious way. It is computed that $| V(G\times H)|=|V(G)|.|V(H)|$ and $|E(G\times H)|=2|E(G)|.|E(H)|$.  We recall that graphs with no pairs of vertices with the same open neighborhoods are called $R$-thin.   In continue, we need the following theorem:

\begin{theorem}{\rm \cite{Sandi}}\label{autKronecker} Suppose $\varphi$ is an automorphism of a connected non-bipartite $R$-thin graph
	$G$ that has a prime factorization $G = G_1 \times G_2 \times \ldots \times  G_k$. Then there exists a permutation
	$\pi$ of $\{1, 2, \ldots , k\}$, together with isomorphisms $\varphi_i : G_{\pi(i)}\rightarrow G_i$, such that
	$$\varphi(x_1, x_2, \ldots , x_k) = (\varphi_1(x_{\pi(1)}), \varphi_2(x_{\pi(2)}), \ldots , \varphi_k(x_{\pi(k)})).$$
	
\end{theorem}
\section{Distinguishing number of Kronecker product of two graphs}
We begin with the distinguishing number of Kronecker product of complete graphs.
\begin{theorem}
	Let $k, n$, and $d$ be integers so that $d\geqslant 2$ and $(d-1)^k < n \leqslant d^k$. Then
	\begin{equation*}
	D(K_k \times K_n) =\left\{
	\begin{array}{ll}
	d & \textsl{if $n\leqslant d^k - \lceil {\rm log}_dk  \rceil -1$}\\
	d + 1& \textsl{if $n\geqslant d^k - \lceil {\rm log}_dk  \rceil +1$}
	\end{array}\right.
	\end{equation*}
	
	If $n=d^k - \lceil {\rm log_dk}  \rceil$, then $D(K_k \times K_n) $ is either $d$ or $d + 1$.
\end{theorem}
\proof It is easy to see that   $K_k \times K_n$ is the complement of Cartesian product $K_k \Box K_n$. By Theorem 1.1 in \cite{Imrich},   $D(K_k\Box K_n)=\left\{
\begin{array}{ll}
d & \textsl{if $n\leqslant d^k - \lceil {\rm log}_dk  \rceil -1$}\\
d + 1& \textsl{if $n\geqslant d^k - \lceil {\rm log}_dk  \rceil +1$}
\end{array}\right.$, so we have the result. \qed

It is known that connected non-bipartite graphs have unique prime factor decomposition with respect to the Kronecker product \cite{McKenzie}. If such a graph $G$ has no pairs $u$ and $v$ of vertices with the same open neighborhoods, then the structure of automorphism group of $G$ depends on that of its prime factors exactly as in the case of the Cartesian product. As said before graphs with no pairs of vertices with the same open neighborhoods are called $R$-thin and it can be shown that a Kronecker product is $R$-thin if and only if each factor is $R$-thin.

\begin{theorem}
	Let $G$ and $H$ be two simple connected, relatively prime graphs, non-bipartite $R$-thin graphs, then $D(G\times H)= D(G\Box H)$.
\end{theorem} 
\proof  By hypotheses and Theorems \ref{autoCartesian} and \ref{autKronecker}, it can be concluded that ${\rm Aut}(G\times H)={\rm Aut}(G\Box H)$. Therefore $D(G\times H)= D(G\Box H)$.\qed

Imrich and  Kla\v{v}zar in \cite{Imrich and  Klavzar} proved that the distinguishing number of  $k$-th power  with respect to the Kronecker product of a non-bipartite, connected, $R$-thin graph different from $K_3$ is two.

\begin{theorem}{\rm \cite{Imrich and  Klavzar}}
	Let $ G$ be a nonbipartite, connected, $R$-thin graph different from $K_3$
	and $\times G^k$ the $k$-th power of $G$ with respect to the Kronecker product. Then $D(\times G^k) = 2$ for $k \geqslant 2$. For the case $G = K_3$ we have $D(K_3 \times K_3) = 3$ and $D(\times K_3^k) = 2$ for $k \geqslant 3$.
\end{theorem}

Now we want to obtain the distinguishing number of  Kronecker product of two complete bipartite graphs. We need the following lemma: 

\begin{lemma} {\rm\cite{Jha}} \label{2.2} 
	If $G = (V_0\cup V_1,E)$ and $H = (W_0\cup W_1,F)$ are bipartite
	graphs, then $(V_0\times  W_0)\cup (V_1 \times W_1)$ and $(V_0 \times W_1)\cup (V_1 \times W_0)$ are vertex sets of the two components of $G \times H$.
\end{lemma} 
\begin{proposition}\label{complete bipartite}
	If $K_{m,n}$ and $K_{p,q}$ are complete bipartite graphs such that $q\geqslant p$ and $m\geqslant n$ then the distinguishing number of $K_{m,n}\times K_{p,q}$ is 
	\begin{equation*}
	D(K_{m,n}\times K_{p,q})=\left\{
	\begin{array}{ll}
	mq+1&  m=n, p=q\\
	mq& {\rm otherwise.} 
	\end{array}\right.
	\end{equation*}
\end{proposition}
\proof The Kronecker product $K_{m,n}\times K_{p,q}$ is disjoint union of two complete bipartite graphs $K_{mp,nq}$ and $K_{mq,np}$ by Lemma \ref{2.2}. Hence if $m\neq n$ and $p\neq q$, then $K_{mp,nq}$ and $K_{mq,np}$ are the two non-isomorphic graphs, and so $D(K_{m,n}\times K_{p,q})={\rm max}\{D(K_{mp,nq}), D(K_{mq,np})\}=mq$. If $m=n$ or $p=q$, then $K_{mp,nq}$ and $K_{mq,np}$ are isomorphic to $K_{mp,mq}$ or $K_{mp,np}$, respectively. In fact
\begin{equation*}
K_{m,n}\times K_{p,q}=\left\{
\begin{array}{ll}
K_{mp,mq}\cup K_{mp,mq}&  m=n, p\neq q \\
K_{mq,nq}\cup K_{mq,nq}&  m\neq n, p =q\\
K_{mq,mq}\cup K_{mq,mq}&  m=n, p =q.
\end{array}\right.
\end{equation*}

Thus using the value of the distinguishing number of complete bipartite graphs we have
\begin{equation*}
D(K_{m,n}\times K_{p,q})=\left\{
\begin{array}{ll}
mq&  m=n, p\neq q~or~ m\neq n, p =q\\
mq+1&  m= n, p =q.
\end{array}\right.
\end{equation*}

Therefore the result follows.\qed
\begin{corollary}
	Let $m,n\geqslant 3$ be two integers. The distinguishing number of Kronecker product of star graphs $K_{1,n}$ and $K_{1,m}$ is $D(K_{1,n}\times K_{1,m})=mn$.
\end{corollary}

The following result shows that the distinguishing number of Kronecker product of complete bipartite graphs is an upper bounds for the distinguishing number of Kronecker product of bipartite graphs.
\begin{corollary}
	If $G=(V_0\cup V_1,E)$ and $H=(W_0\cup W_1,E)$ are bipartite graphs such that $|V_0|=m$, $|V_1|=n$, $|W_0|=p$, and $|W_1|=q$, then 
	$D(G\times H)\leqslant D(K_{m,n}\times K_{p,q})$.
\end{corollary}
\proof It is sufficient to note that ${\rm Aut} (G\times H)\subseteq {\rm Aut}(K_{m,n}\times K_{p,q})$, and $G\times H$ and $K_{m,n}\times K_{p,q}$ have the same size. Now we have the result by Proposition \ref{complete bipartite}.\qed

Before we prove the next result we need some additional information  about the distinguishing number of  complete multipartite graphs. Let $K_{{a_1}^{j_1},{a_2}^{j_2},\ldots, {a_r}^{j_r}}$  denote the complete multipartite graph that has $j_i$ partite sets of size $a_i$ for $i = 1, 2,\ldots , r$ and $a_1 > a_2 > \ldots > a_r$. 
\begin{theorem}\label{Theorem 2.4 of Collins and A. N. Trenk}{\rm \cite{Collins and A. N. Trenk}}
	Let $K_{{a_1}^{j_1},{a_2}^{j_2},\ldots, {a_r}^{j_r}}$ denote the complete multipartite graph that has $j_i$ partite sets of size $a_i$ for $i = 1, 2,\ldots , r$, and $a_1 > a_2 > \ldots > a_r$. Then 
	\begin{equation*}
	D(K_{{a_1}^{j_1},{a_2}^{j_2},\ldots, {a_r}^{j_r}}) ={\rm min}\{p : {p \choose a_i} \geqslant j_i ~{\rm for~all~i}\}.
	\end{equation*}
\end{theorem}

\begin{theorem}
	If $G$ and $H$ are two simple connected, relatively prime  graphs such that $G\times H$ has $j_i$ $R$-equivalence classes of sie $a_i$ for $i=1,\ldots , r$, and $a_1 > a_2 >\ldots > a_r$ then
	\begin{equation*}
	D(G\Box H)\leqslant D(G\times H) \leqslant {\rm min}\{p: {p\choose a_i}\geqslant j_i ~{\rm for~ all~ i}\}.
	\end{equation*}
\end{theorem}
\proof  Since ${\rm Aut}(G\Box H)\subseteq {\rm Aut}(G\times H)$, so $D(G\Box H)\subseteq D(G\times H)$. To prove the second inequality, it is sufficient to consider each $R$-equivalence classes of $G\times H$ as a partite set. Thus graph $G\times H$ can be considered as multipartite graph that has $j_i$ partite sets of size $a_i$ such that every two partite sets of this multipartite graph is complete bipartite or there exists no edge between the two partite sets. So the automorphism group of  this multipartite graph is subset of the automorphism group of complete multipartite graph with the same partite sets. Therefore $D(G\times H) \leqslant D(K_{a_1^{j_1}, \ldots , a_r^{j_r}})$, and the result follows by Theorem \ref{Theorem 2.4 of Collins and A. N. Trenk}.\qed

By using the concept of the Cartesian skeleton we can obtain an upper bound for Kronecker product of $R$-thin graphs. For this purpose we need the following preliminaries frome \cite{Sandi}: The Boolean square of a graph $G$ is the graph $G^s$ with $V(G^s) = V (G)$ and $E(G^s) =
\{xy ~|~ N_G(x)\cap N_G(y) \neq \emptyset\}$. An edge $xy$ of the Boolean square $G^s$ is dispensable if it is a loop, or if
there exists some $z \in V(G)$ for which both of the following statements hold:
\begin{enumerate}
	\item[(1)] $ N_G(x) \cap N_G(y) \subset N_G(x) \cap N_G(z) ~{\rm or}~ N_G(x) \subset N_G(z) \subset N_G(y)$.
	\item[(2)] $N_G(y) \cap N_G(x) \subset N_G(y) \cap N_G(z) ~{\rm or}~ N_G(y) \subset N_G(z)\subset N_G(x)$.
\end{enumerate}

The Cartesian skeleton $S(G)$ of a graph $G$ is the spanning subgraph of the
Boolean square $G^s$ obtained by removing all dispensable edges from $G^s$.
\begin{proposition}{\rm \cite{Sandi}}\label{skeleton} If $H$ and $K$ are $R$-thin graphs without isolated vertices, then $S(H \times K) =S(H) \Box S(K)$.
\end{proposition}

\begin{proposition}{\rm \cite{Sandi}}\label{isomorphskeleton}  Any isomorphism $\varphi : G \rightarrow H$, as a map $V(G) \rightarrow V(H)$, is also an
	isomorphism $\varphi : S(G) \rightarrow S(H)$.
\end{proposition}

Now we are ready to give an upper bound for Kronecker product of $R$-thin graphs.
\begin{theorem}
	If $G$ and $H$ are $R$-thin graphs without isolated vertices, then $D(G\times H)\leqslant D(S(G) \Box S(H))$. 
\end{theorem}
\proof By Proposition \ref{isomorphskeleton} we have ${\rm Aut}(G\times H)\subseteq {\rm Aut}(S(G\times H))$, and so $D(G\times H)\subseteq D(S(G\times H))$. The result obtains immediately from Proposition \ref{skeleton}. \qed

\section{Distinguishing index of Kronecker product of two graphs}
In this section we investigate the distinguishing index of Kronecker product of two graphs. Let us start with the Kronecker product power of $K_2$. It can be seen that $\times K_2^n$ is disjoint union of $2^{n-1}$ number of $K_2$, and so $D'(\times K_2^n)=2^{n-1}$ for $n\geqslant 2$. 

\begin{figure}
	\begin{center}
		\includegraphics[width=0.6\textwidth]{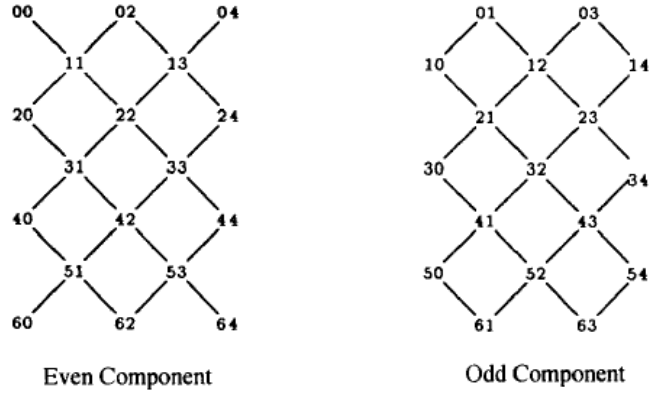}
		\caption{\label{fig1}Graph $P_7 \times P_5$.}
	\end{center}
\end{figure}

Let $ij$ be the notation  of the vertices the Kronecker product $P_m\times P_n$  of two paths of order $m$ and $n$ where $0\leqslant i \leqslant m-1$ and $0\leqslant j \leqslant n-1$. From \cite{Jha2}, we know that  
$P_m\times P_n$ is bipartite, the number of vertices in even component  of $P_m\times P_n$ is $\lceil mn/2  \rceil$ while that in odd component is $\lfloor mn/2  \rfloor$ (see Figure \ref{fig1}). 
It can be easily computed that distinguishing index of $P_m\times P_n$ is two, unless $D'(P_3\times P_2)=3$ and $D'(P_3\times P_3)=4$, because $P_3\times P_2$ is disjoint union $P_3\cup P_3$, and $P_3\times P_3$ is disjoint union $K_{1,4}\cup C_4$.

The distinguishing index of the square of cycles and for arbitrary power of complete graphs with respect to the Cartesian, Kronecker and strong  product has been considered by Pil\'sniak \cite{pilsniak}. In particular, she proved that $D'(\times C_m^2)=2$ for the odd value of $m\geqslant 5$, and $D'(\times K_n^r)=2$ for any $n\geqslant 3$ and $r\geqslant 2$. 

Let us state and prove  the following lemma concerning the Kronecker product $K_2\times H$.

\begin{lemma}\label{K_2H}
	If $H$ is a graph with $D'(H)=d$, then $D'(K_2\times H)\leqslant d+1$. If $H$ is bipartite then  $d\leqslant D'(K_2\times H)\leqslant d+1$.
\end{lemma}
\proof Let $f$ be an automorphism of bipartite graph $K_2\times H$ with partite sets $\{(v_1,x) | x\in V(H)\}$ and $\{(v_2,x) | x\in V(H)\}$. Since $f$ preserves the adjacency and non-adjacency relations, so $f(v_i,x)=(v_i,\varphi (x))$ for all $x\in V(H)$ or $f(v_i,x)=(v_j,\varphi (x))$ for all $x\in V(H)$ where $i,j\in \{1,2\}$, $i\neq j$, and $\varphi \in {\rm Aut}(H)$.

Let $L$ be a distinguishing edge labeling of $H$. If $(v_1,h_1)(v_2,h_2)$ be an arbitrary edge of $K_2\times H$, then we assign it the label of the edge $h_1h_2$ in $H$. Now suppose that $hh'$ is an edge of $H$ and fix it. We change the label of the edge $(v_1,h)(v_2,h')$ to a new label. If $f$ is an automorphism of $K_2\times H$ preserving the labeling, then with respect to the label of the two edges $(v_1,h)(v_2,h')$ and $(v_1,h')(v_2,h)$ we have $f(v_i,h)=(v_i,\varphi (x))$ for $i=1,2$ and some $\varphi \in {\rm Aut}(H)$. On the other hand $\varphi$ is the identity, because we labeled the edges of $K_2\times H$ by the distinguishing edge labeling $L$ of $H$. Therefore this labeling is distinguishing.
If $H$ is bipartite then $K_2\times H = H\cup H$, and so the result follows.\qed
\begin{remark}
Let $(G, \phi)$ denote the labeled version of $G$ under the labeling $\phi$. Given two distinguishing $k$-labelings
$\phi$ and $\phi'$ of $G$, we say that $\phi$ and $\phi'$ are equivalent if there is some automorphism of $G$ that maps $(G, \phi)$
to $(G, \phi')$. We denote by $D(G, k)$ the number of inequivalent $k$-distinguishing labelings of $G$  which was first considered by Arvind and Devanur \cite{Arvind and Devanur} and Cheng \cite{ChengD(Gk)} to determine
the distinguishing numbers of trees.
 In Lemma \ref{K_2H}, if $H$ is bipartite and $D(H,d)=1$ then  $ D'(K_2\times H)=d+1$; otherwise, i.e., if $H$ is bipartite and $D(H,d)>1$ then  $ D'(K_2\times H)= d$.
\end{remark}
\begin{proposition}
	If $m\geqslant 4$ and $n\geqslant 2$, then $D'(P_m\times K_{1,n})=n$. Also $D'(P_2\times K_{1,n})=n+1$ and $D'(P_3\times K_{1,n})=2n$.
\end{proposition}
\proof
The consecutive vertices of $P_m$ are denoted by $0,1,\ldots , m-1$ and the central vertex of $K_{1,n}$ by $0$, and its pendant vertices by $1,\ldots ,n$. Since $K_{1,n}$ and $P_m$ are connected and bipartite, their Kronecker product consists of two connected components (see Figure \ref{fig2,3}).
\begin{figure}[ht]
	\begin{center}
		\begin{minipage}{6.3cm}
			\includegraphics[width=1.00\textwidth]{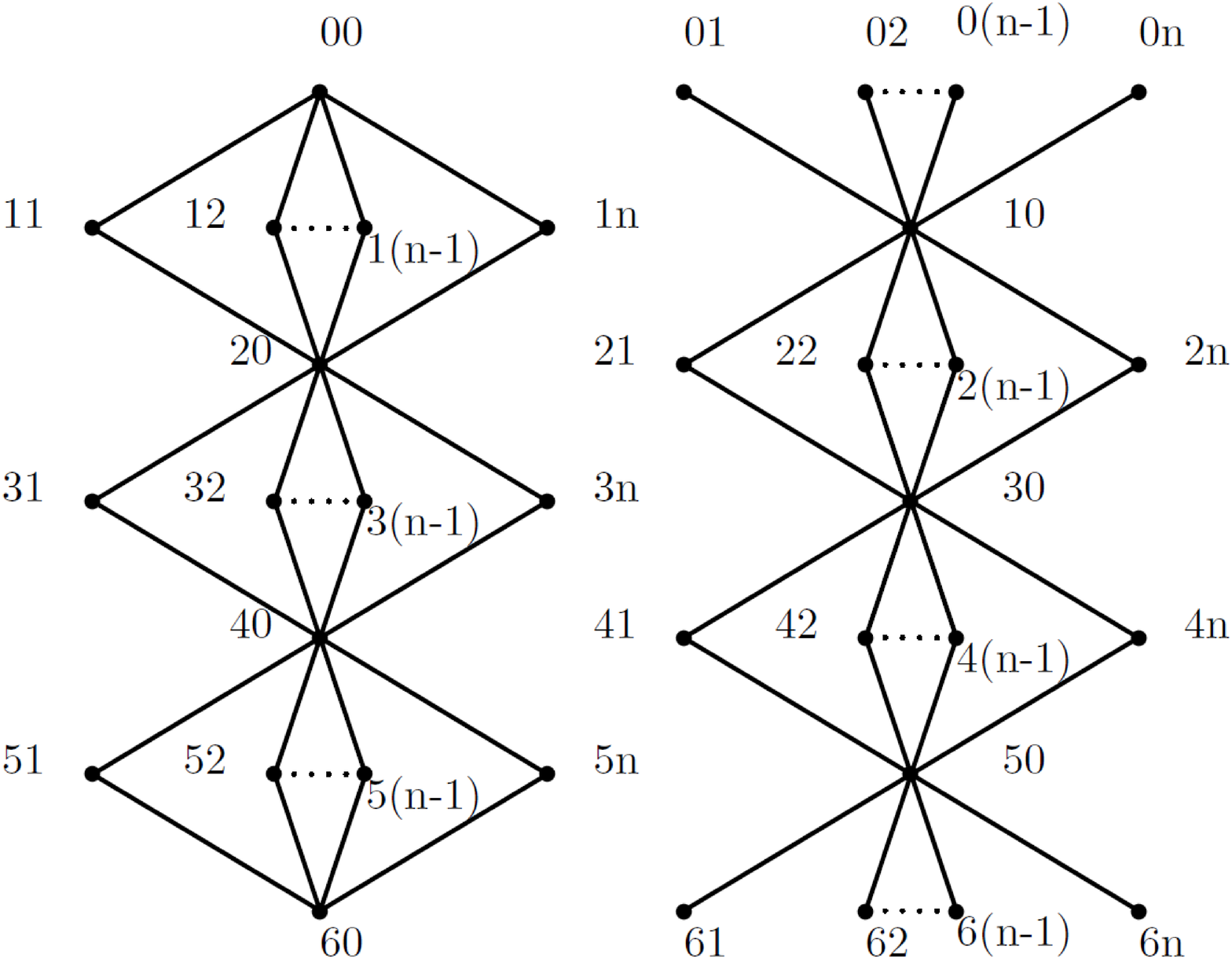}
		\end{minipage}
		\hspace{1cm}
		\begin{minipage}{6.3cm}
			\includegraphics[width=1.00\textwidth]{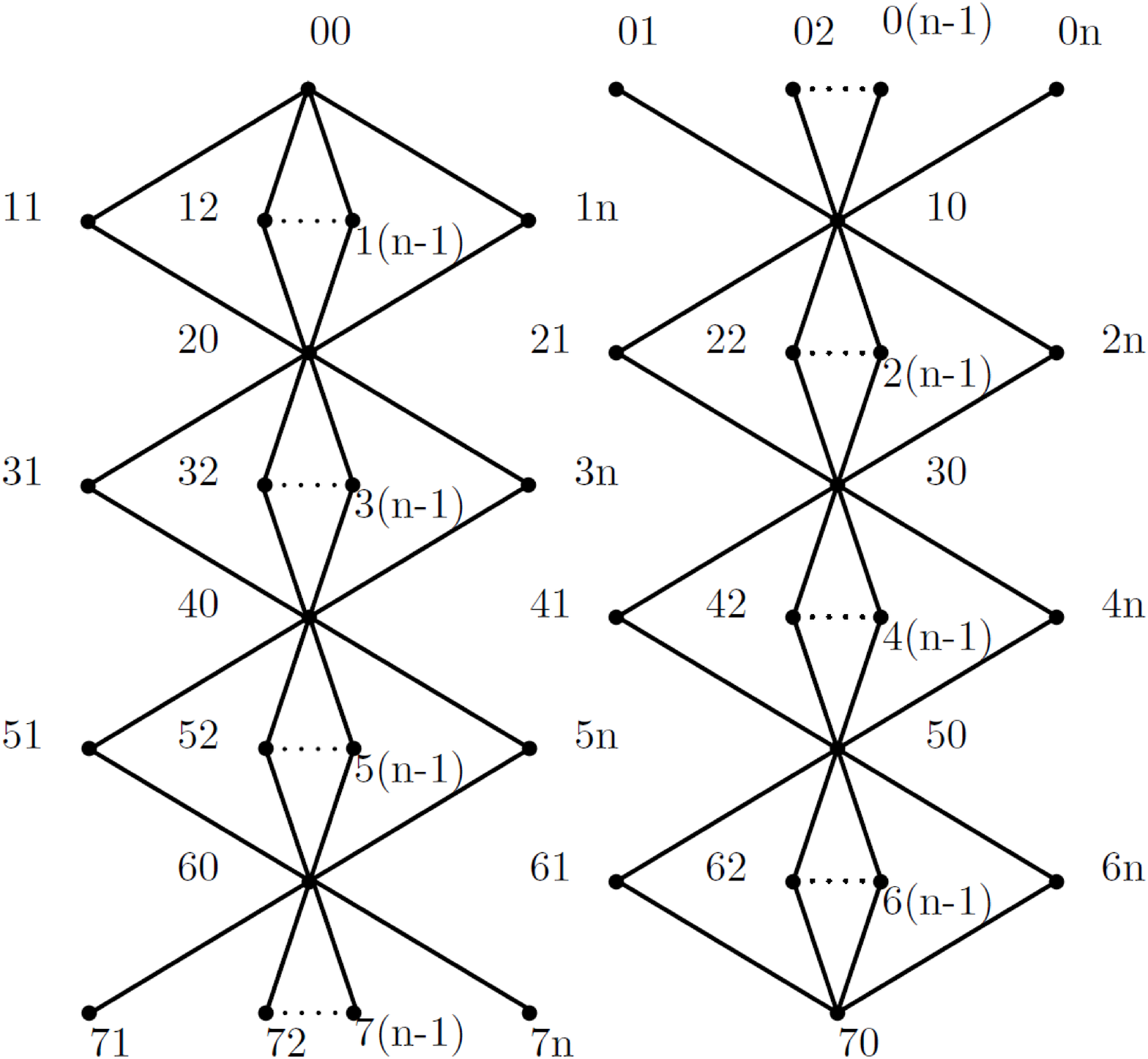}
		\end{minipage}
		\caption{\label{fig2,3}Graphs $K_{1,n}\times P_7$ and $K_{1,n}\times P_8$. }
	\end{center}
\end{figure}

If $m$ is odd, then we label the edges $(10,0i)$ and $(j0,(j+1)i)$ with label $i$ where $1\leqslant i \leqslant n$ and $0\leqslant j \leqslant m-1$. We label the remaining edges with an arbitrary label, say 1. Using Figure  \ref{fig2,3} and regarding to the number of pendant vertices we can obtain that this labeling is distinguishing and $D'(K_{1,n}\times P_m)=n$. If $m$ is even, then we label the edges $(10,0i)$, $(00,1i)$, and $(j0,(j+1)i)$ with label $i$ for $1\leqslant i \leqslant n$ and $2\leqslant j \leqslant m-2$. Also we label the edges $((m-2)0,(m-3)i)$ with label 1 and  the edges $(10,2i)$ with label 2 for $1\leqslant i \leqslant n$. We label the remaining edges with an arbitrary label, say 1. With a similar argument  we can conclude that this labeling is distinguishing and $D'(P_m\times K_{1,n})=n$. 
Since $K_{1,n}\times P_2$ is disjoint union $K_{1,n}\cup K_{1,n}$,  and $K_{1,n}\times P_3$ is disjoint union $K_{1,2n}\cup G$ where $G$ is a graph with $D'(G) \leqslant n$, so $D'(K_{1,n}\times P_2)=n+1$ and $D'(K_{1,n}\times P_3)=2n$. \qed

\begin{proposition}
	If $n\geqslant m\geqslant 3$, then $D'(K_{1,n}\times K_{1,m})=nm$.
\end{proposition}
\proof By definition, $K_{1,n}\times K_{1,m}$ is disjoint union $K_{1,nm}\cup K_{n,m}$, and hence $D'(K_{1,n}\times K_{1,m})={\rm max}\{D'(K_{1,nm}),D'(K_{n,m})\}$. On the other hand $D'(K_{n,m})\leqslant \lceil \sqrt[m]{n} \rceil +1$ (by Corollary 20 from \cite{nord}). Therefore $D'(K_{1,n}\times K_{1,m})=D'(K_{1,nm})=nm$.\qed

\begin{figure}
	\begin{flushleft}
		
		\includegraphics[width=1.0\textwidth]{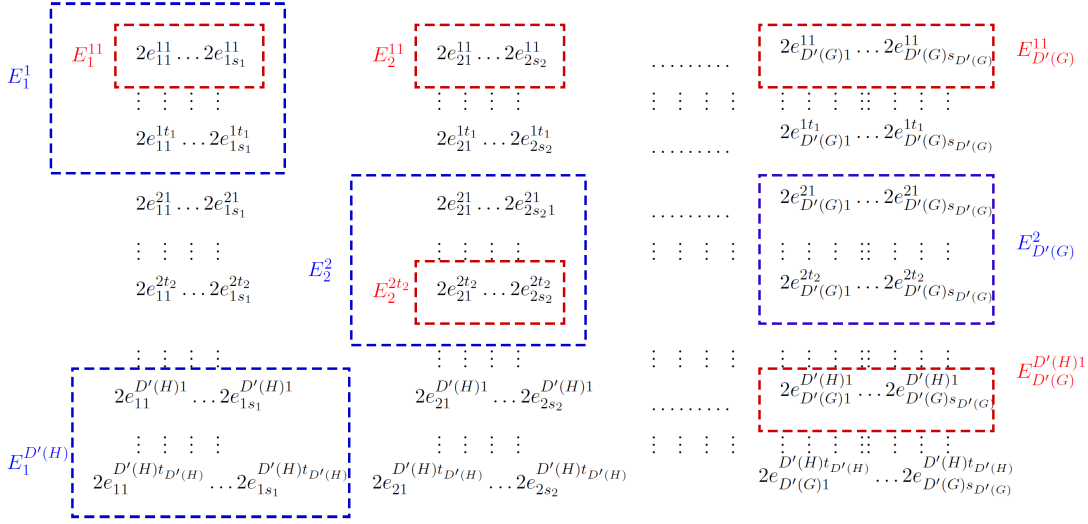}
		\caption{\label{fig4}The partition of the edge set of $X$ in Theorem \ref{thmD'}.}
	\end{flushleft}
\end{figure}

\begin{theorem}\label{thmD'}
	Let $X$ be a connected non-bipartite $R$-thin graph which has a prime factorization $X=G\times H$ where $G$ and $H$ are simple and ${\rm max}\{D'(G),D'(H)\}\geqslant 2$. Then $D'(X)\leqslant D'(K_{D'(G),D'(H)})$.
\end{theorem}
\proof Let the sets $\{a_{i1},\ldots ,a_{is_i}\}$ where $1\leqslant i \leqslant D'(G)$ and the sets $\{b_{j1},\ldots ,b_{jt_j}\}$ where $1\leqslant j \leqslant D'(H)$ be the partitions of the edges set $G$ and $H$ by its distinguishing edge labeling, respectively, i.e., the label of $a_{ij}$ is $i$ for $1\leqslant j\leqslant s_i $ and  the label of $b_{ij}$ is $i$ for $1\leqslant j\leqslant t_i$. If $a_{ij}$ is the edge of $G$ between $v$ and $v'$, and $b_{pq}$ is the edge of $H$ between $w$ and $w'$, then $2e_{pq}^{ij}$ means the two edge $(v,w)(v',w')$ and $(v,w')(v',w)$ of $X$. Regarding to above partitions we can partition the edge set of $X$ as the sets $E_p^{ij}=\{2e_{pq}^{ij}~|~ 1 \leqslant q \leqslant t_p\}$, where $1\leqslant i \leqslant D'(G)$, $1\leqslant j\leqslant s_i $, and $1\leqslant p \leqslant D'(H)$. Set $E_p^i =\bigcup_{j=1}^{t_i} E_p^{ij}$ (see Figure \ref{fig4}).
If we assign the all elements of $E_p^i$, the label of the edge $z_iz'_p$ of complete bipartite graph $K_{D'(G),D'(H)}$ with parts $\{z_1,\ldots , z_{D'(G)}\}$ and $\{z'_1,\ldots , z'_{D'(H)}\}$, then we have a distinguishing edge labeling of $X$ by Theorem \ref{autKronecker}, and so the result follows.\qed

\begin{theorem}\label{D'=1}
	Let $X$ be a connected non-bipartite $R$-thin graph which has a prime factorization $X=G\times H$ where $G$ and $H$ are simple and $D'(G)=D'(H)=1$. Then $D'(X)\leqslant 2$.
\end{theorem}
\proof  If  $G$ and $H$ are non-isomorphic then $|{\rm Aut}(X)|=1$ by Theorem \ref{autKronecker}, and so $D'(X)=1$. Otherwise, $D'(G\times G)=2$ because, $G\times G$ is a symmetric graph and we have a 2-distinguishing edge labeling of it as follows: according to notations of Theorem \ref{thmD'}, if we label the all elements of $E_1^{11}$ with label 1, and all elements of $E_1^{1j}$ with label 2 for every $2\leqslant j \leqslant |E(G)|$, then we have a 2-distinguishing edge labeling of $G\times G$.\qed

\begin{corollary}
If $G$ is  a connected non-bipartite $R$-thin graph with $D'(G)=1$, then $D'(\times G^k)=2$ for any $k\geqslant 2$.
\end{corollary} 
\proof It follows immediately by   induction on $k$ and Theorems \ref{thmD'} and \ref{D'=1}.\qed

\end{document}